\magnification=\magstep1
\input amstex
\documentstyle{amsppt}
\pagewidth{16 true cm}
\pageheight{23 true cm}
\topmatter
\TagsOnRight
\NoBlackBoxes
\rightheadtext { }
\title
Igusa's local zeta functions of semiquasihomogeneous polynomials
\endtitle
\rightheadtext{ }
\leftheadtext{}
\author
 W. A. Z\'u\~niga-Galindo*
\endauthor
\address
\noindent Universidad Aut\'{o}noma de Bucaramanga,
Laboratorio de Computo Especializado,
A.A. 1642, Bucaramanga,
Colombia
\endaddress
\email \noindent wzuniga\@bumanga.unab.edu.co
\endemail
\thanks
* This work was supported by  COLCIENCIAS, contract \# 063-98
\endthanks
\abstract
In this  paper, we  prove the  rationality of Igusa's local zeta functions of semiquasihomogeneous polynomials with  coefficients in a non-archimedean local field $K$.  The proof of this result is based on Igusa's stationary phase formula and some ideas on N\'{e}ron $\pi-$desingularization.

\endabstract
\endtopmatter

\document
\openup 3pt

\head
 {\bf 1. Introduction}
\endhead

 Let $K$ be a  non-archimedean  local field, and let $\Cal{O}_{K}$ be  the   ring of integers 
of $K$. Let $\pi$ be  a uniformizing parameter of $K$, and let the residue field of $K$ be 
   $\Bbb{F}_q$ a finite field  with $q = p^{r}$ elements.  Let $v$ denote the valuation of $K$ such that $v(\pi) =1$. For $x \in K$, let   $ |x|_{K} = q^{-v(x)}$. Let $\mid dx \mid$ be the  Haar measure on $K^{n}$ so normalized that the measure of
$\Cal{O}^{n}_{K}$ is equal to one.
Let  $f(x) \in  K[x]$, $x =(x_1,..,x_n)$, the  Igusa local zeta function associated to  $f$ is 
defined by

$$ Z(f,s)= \int_{\Cal{O}^{n}_{K}}|f(x)|^{s}_K \mid dx \mid , \,\,\,\,\,\,\,\,  s \in \Bbb{C}, \,\, Re(s) >0.$$

\noindent The local zeta function  $Z(f,s)$ is a holomorphic
 function  on the semiplane $Re(s) >0$. In  the case
of  $K$ having  characteristic zero, Igusa ([7], [8]) and Denef
 ([3]) proved that  $Z(f,s)$ is a rational function of $q^{-s}$.
At the present time, the techniques used by Igusa (resolution of singularities)  and
 Denef (elimination of quantifiers in  $\Bbb{Q}_p$)
 are not available in positive characteristic,
so in this case the rationality of $Z(f,s)$ is still an open problem.

\noindent The local zeta function contains information about the number of solutions of the congruence
 $f(x) \equiv 0 $ mod $\pi^{j} \Cal{O}_K$, (see e.g. [4]). More precisely, if 
$$ N_j := \text{Card}\{ x \in \left(\Cal{O}_K / \pi^{j}\Cal{O}_K\right)^{n} \,\ \mid  f(x) \equiv 0 \,\ \text{mod} \,\
  \pi^{j}\Cal{O}_K  \},$$

\noindent and $P(t )$
 is the  Poincar\'{e} series 
 $ \,\,\, P(t )= \sum_{j \geqq 0}^{\infty} N_j (q^{-n}t )^{j}$, then 

  $$Z(f,s) = P(q^{-s})-q^{s}(P(q^{-s})-1).$$

 In this paper, we shall study the local zeta functions of  semiquasihomogeneous polynomials with an absolutely algebraically isolated   singularity at the origin of $K^{n}$.  

\noindent Let $f(x)$   be a polynomial with coefficients in $K$, and $V_f$ the corresponding $K-$hy\-persur\-face. We call a point $P\in K^{n}$  an {\it 
 absolutely algebraically isolated singularity }of $V_f(K)$, if the only solution of  the system of equations
$$
 f(x) = 
\frac{\partial f}{\partial x_1}(x)= ....= \frac{\partial f}{\partial x_n}(x)= 0,
$$
\noindent over  an algebraic fixed closure of $K$, is the point $P$. 

 \noindent  Let $\alpha_1,..,\alpha_n$ be $n$
  relatively prime and positive integers. A polynomial $f(x) \in K[x]$ is called a {\it quasihomogeneous polynomial}
  of weight  $d$ and exponents $ \alpha_1,..,\alpha_n$, if it satisfies:

$$f(t^{\alpha_1}x_1,..,t^{\alpha_n}x_n) = t^{d}f(x),\,\,\,\,\,\,\, \text{for every} \,\,  t \in 
K,$$ 

\noindent and the origin of $K^{n}$ is an absolutely algebraically isolated
singularity of  $K-$hyper\-surface
$V_f$.

A polynomial $F(x)$ is called a {\it semiquasihomogeneous polynomial} if it has the  form $f(x)+\sum b_i e_i(x) \in K[x]$, where $f(x)$ is a quasihomogeneous
polynomial, and each monomial $e_i(x)=x_1^{m_1}...x_n^{m_n}$ satisfies  $\sum \alpha_im_i >d$, and 
 the origin of $K^{n}$ is an absolutely algebraically isolated
singularity of  $K-$hypersurface
$V_F$. We call the polynomial  $f(x)$ the quasihomogeneous part of $F(x)$.

  We  put 
$|\alpha|= \sum_{i}\alpha_i$, for any  $\alpha=( \alpha_1,..,\alpha_n)
 \in \Bbb{Z}^{n}$.
We use the notation $Z(f,D,s)$ for the integral $\int_{D}|f(x)|^{s}_K \mid dx \mid $.
In the case of  $D=\Cal{O}^{n}_{K}$, we  use the simplified notation $Z(f,s)$.

The main result of this paper is the following:

\proclaim{Theorem 3.5}
  Let $F(x) \in K[x]$ be a semiquasihomogeneous
 polynomial whose  quasihomogeneous part $f(x)$ has 
 weight  $d$
and exponents $ \alpha_1,..,\alpha_n$.
 Then   Igusa's local zeta function of $F(x)$ is a rational function of $q^{-s}$. More 
precisely,

$$Z(F,s) = \frac{L(q^{-s})}{ (1-q^{-1}q^{-s})(1-q^{-|\alpha|}q^{-ds})}
,\tag{1.1}$$

\noindent where $\alpha= (\alpha_1,...,\alpha_n)$.  Furthermore, the 
polynomial $L(q^{-s})$ can be computed effectively.
\endproclaim

 If in addition the polynomial $F$ is nondegenerate  for its Newton diagram  and if $K$ has 
characteristic zero, then a very different way of calculating $Z(F,s)$ is given in [5].
The  proof  of theorem 3.5  gives an effective method to compute the local
  zeta functions of semiquasihomogeneous polynomials. 

\noindent We say that   a singular point $ \bar{P} \in V_{\bar{f}}( \Bbb{F}_q)$, where $\overline{f }$ is the reduction modulo $\pi$ of $f$, is a {  non-liftable singularity} of  the hypersurface $V_f$, if
for every singular point $Q \in V_f(K)$, $Q \in \Cal{O}_K^{n}$, the reduction modulo $\pi$ of $Q$ is different of $\overline{P}$.
 The proof of  theorem 3.5 shows that the numerator of the zeta
   function $Z(F,s)$ depends on the non-liftable singularites of the closed fiber of the
    hypersurface $V_F$, and the denominator depends on the
     singularity of the generic fiber of $V_F$. More precisely,
      the denominator
     depends on Newton's diagram of $F(x)$.  
  In the  proof of theorem 3.5, we use Igusa's formula of stationary phase
   for $\pi$-adic integrals ([10]) and some  ideas  on
    N\'{e}ron $\pi$-desingularization  ( [13], sect. 17, 18).

As a consequence of theorem 3.5, we obtain the following three corollaries.

\proclaim{Corollary 3.6}
\noindent Let $K$ be a global field and  $F(x) \in K[x]$ be a semiquasihomogeneous
 polynomial whose quasihomogeneous part $f(x)$  has  weight $d$
and exponents $ \alpha_1,..,\alpha_n$. 
Then for  every non-archimedean valuation $v$ of $K$,  Igusa's local zeta function of $F(x)$ on the 
completion $K_v$ of $K$ is a rational function of form (1.1). If $K$ is a number field and  $F(x)$ is non-degenerate for
 its Newton's diagram, then  the real parts of the poles of the zeta function $Z(F,s)$ are roots of
  the Bernstein polynomial of $F(x)$.
\endproclaim

For the definition of the Bernstein polynomial
 and its computation in the non-degenerate case see reference [2]. The last part of corollary
 3.6 is a special case of a more
   general result due to     Loeser (cf. [11], thm. 5.5.1).

\proclaim{Corollary 3.7}
\noindent Let $K$ be a global field and  $\Cal{O}_K$ its ring of integers. Let $F(x) \in \Cal{O}_K[x]$ be a semiquasihomogeneous polynomial whose quasihomogeneous part $f(x)$ has weight $d$ and exponents $\alpha_1,...,\alpha_n$. Then for every non-archimedean 
valuation $v$ of $K$, the number of solutions $N_j(F,v)$ of the congruence
$$ F(x)
 \equiv 
0 
\,\ 
\text{mod} 
\,\  \pi^{j}\Cal{O}_{K_{v}},$$
where  $\Cal{O}_{K_{v}}$  is the ring of integers of the completion $K_v$, satisfies 
$$
\text{l\'\i m sup}_{j\to \infty}N_{j}(F,v)^{1/j} \le 
\cases
 q^{n- \mid \alpha \mid  / d }\,\,\,\,\ \text{if}  \,\,\ \mid \alpha \mid  / d \le 1, \,\,\ \\
q^{n- 1 } \,\,\,\,\,\,\,\,\,\,\,\  \text{if} \,\,\,\, \mid \alpha \mid  / d > 1.\\
\endcases
$$

\endproclaim

 The zeta functions of plane curves, with only an absolutely analytically irreducible
 singularity at the origin, have  been extensively studied when the characteristic of $K$ is zero,
  by Igusa ([9]), Meuser ([12]), among others. Let $f(x,y) \in K[x,y]$ be
  an absolutely analytically irreducible polynomial. Thus there exist $(\alpha_1,
  \alpha_2) \in \Bbb{N}^{2}$, relatively prime integers, and an integer $d$, such that $f(x,y) = f_d(x,y)
  +g(x,y)$ and the monomials $x^{n} y^{m}$ of $f_d(x,y)$ and $g(x,y)$ satisfy
  $\alpha_1 n +\alpha_2 m =d$ and $\alpha_1 n +\alpha_2 m >d$, respectively.
  The origin of $K^{2}$ is an absolutely algebraically isolated singularity 
of the
   plane curve $V_f$. If it is also valid for   $V_{f_d}$, then $f$ is a semiquasihomogeneous polynomial   in the sense of the
  definition given for us. As a consequence of theorem 3.5, we  obtain a precise
   description of the local zeta function associated with this type of polynomials.

   \medskip
\proclaim{Corollary 3.8}

\noindent Let   $f(x,y) \in K[x,y]$ be an absolutely  analytically irreducible polynomial,
such that the origin of $K^{2}$ is an absolutely algebraically isolated singularity of
 the plane
curve  $V_{f_d}$.  Then  the Igusa local zeta function $Z(f,s)$ is a rational function of
$q^{-s}$ of  form (1.1).
\endproclaim

\head
{\bf 2. Preliminaires}
\endhead

\noindent In [10] Igusa introduced the stationary phase formula for $\pi-$adic integrals  and 
suggested that a closer examination of this formula might lead to a proof of the rationality of 
$Z(f,s)$ in any characteristic. The  above suggestion has been our main motivation for
this paper.  In this section we review   Igusa's stationary
phase formula  and some  ideas on N\'{e}ron 
$\pi-$desingularization.

We denote by $\overline{x}$ the image of an element of $\Cal{O}_K$ under the canonical 
homomorphism $ \Cal{O}_K \longrightarrow  \Cal{O}_K / \pi 
\Cal{O}_K \cong \Bbb{F}_q $, i.e, the reduction modulo  $\pi$. Given $f(x) 
 \in \Cal{O}_K[x] $ such that not all its coefficients are in
  $\pi \Cal{O}_K$,  we denote by $\overline{f(x)}$ the polynomial
   obtained by reducing modulo $\pi$  the coefficients of $f(x)$.

For any commutative ring  $A$ and $f(x) \in A[x]$, we denote by  $V_{f}(A)$, the set of $A-$valued points 
of the hypersurface $V_f$  defined by $f$, and  by $Sing_{f}(A)$, the set of $A-$valued singular points of $V_f$, i.e.,
$$Sing_f(A)= \{x \in A ^{n}\mid \, f(x) = 
\frac{\partial f}{\partial x_1}(x)= ....= \frac{\partial f}{\partial x_n}(x)= 0\}.$$

We fix a lifting $R$ of $\Bbb{F}_q$ in $\Cal{O}_K$. Thus, the   set $R$ is mapped bijectively onto
$\Bbb{F}_q$ by the canonical homomorphism $ \Cal{O}_K \longrightarrow  \Cal{O}_K / \pi 
\Cal{O}_K $.

\noindent 

Let $f(x) \in \Cal{O}_K[x]$ be  a polynomial in $n$ variables,
  $ P_1= (y_{1},...,y_{n})\in
  \Cal{O}_K^{n}$,  and  $m_{P_1}=(m_1,...,m_n) \in \Bbb{N}^{n}$. We call a  $K^{n}-$isomorphism
$\psi_{m_{P_1}}(x)$ a {\it dilatation}, if it satisfies $\psi_{m_{P_1}}(x)= (z_1,..,z_n)$, 
$ z_i = y_i+ \pi ^{m_{P_i}}x_i$, for each $i =1,2,..,n$. We define the  {\it  dilatation }of $f(x)$ at $P_1$ induced by $\psi_{m_{P_1}}(x)$, as 

$$f_{P_1}(x):= \pi^{-e_{P_1}}f(\psi_{m_{P_1}}(x)),\tag{2.1}$$

\noindent where   $e_{P_1}$ is the minimun order   of $\pi$ in  the coefficients
of $f(\psi_{m_{P_1}}(x))$.  We call the $K-$hypersurface $V_{f_{P_1}}$ the dilatation of 
$V_f$ at $P_1$ induced by $\psi_{m_{P_1}}(x)$, the number $e_{P_1}$  
{\it  the arithmetic multiplicity  of $f(x)$ at $P_1$} by $\psi_{m_{P_1}}(x)$
  and the set  $ S(f_{P_1})$, the lifting of $Sing_{\overline{f_{P_1}}}(\Bbb{F}_q)$, {\it the first generation of descendants
  of $P_1$}. By dilatating $f_{P_1}$ at each $P_2 \in  S(f_{P_1})$, by   some  $\psi_{m_{P_2}}(x)$, we obtain $\left(f_{P_1 , P_2}, S(f_{P_1 , P_2} ) \right)_{P_2 \in  S(f_{P_1})}$. The union of the sets $\left( S(f_{P_1 , P_2} ) \right)_{P_2 \in  S(f_{P_1})}$ is the second generation of descendants of $P_1$. Given a sequence of dilatations $( \psi_{m_{P_k}}(x)) _{k}$, we define inductively $e_{P_1,...,P_k}$ and $f_{P_1,...,P_k}(x)$, 
$ S(f_{P_1,...,P_k})$  as follows:

$$f_{P_1,...,P_k}(x)= \cases
f(x) \,\,\,\,\ \text{if} \,\,\ k=0,\\
\pi^{-e_{P_1,...,P_k}}f_{P_1,...,P_{k-1}}( \psi _{m_{P_k}} (x)) \,\,\,\,\text{if} \,\,\,\, k\geqq 1,\\
\endcases
\tag{2.2}$$

\noindent where $P_k$ runs through the descendants of $k-1$
generation of $P_1$. The union of the sets $S(f_{P_1,...,P_k})$,
  is called {\it the $k-$generation of descendants of $P_1$}.

 Along this paper we shall use several types of dilatations, i.e., dilatations  with different $m$'s, however the especific value of $m$ will be clear from the context. The dilatations were introduced by N\'{e}ron (cf. [13], sect. 18). These transformations play an important role  in the process of desingularization of the closed fiber of a scheme over  a  discrete valuation ring whose generic fiber is non-singular. A modern exposition of the  N\'{e}ron $\pi -$desingularization can be found in [1], sect. 4.

\noindent Now, we review Igusa's stationary phase formula, from the point of view of the dilatations. For that, we  fix the $m_{P_k}$'s equal to $(1,..,1) \in \Bbb{N}^{n}$ in (2.2).

 Let $\overline{D}$ be a subset of $\Bbb{F}_q^{n}$ and 
  $D$ its preimage under the canonical homomorphism
 $\Cal{O}_K \longrightarrow  \Cal{O}_K / \pi 
\Cal{O}_K \cong \Bbb{F}_q $. Let  $S(f,D)$ denote   the subset of $R^{n}$
 (the set of representatives of
$\Bbb{F}_q^{n}$ in $\Cal{O}_K^{n}$) mapped bijectively to the set
$Sing_{\bar{f}}(\Bbb{F}_q) \bigcap \overline{D}$. We use the simplified notation $S(f)$, in the case of
$D= \Cal{O}_K^{n}$. Also we define:

$$ \nu (\bar{f},D):= q^{-n}\text{Card} \{\overline{P} \in \overline{D} |\,\,  
\overline{P} \notin V_{\bar{f}}(\Bbb{F}_q)\},$$

$$\sigma (\bar{f},D) := q^{-n}\text{Card} \{\overline{P} \in \overline{D} |\,\,  
\overline{P} \,\,\ \text{is a non-singular point of} \,\,\,
 V_{\bar{f}}(\Bbb{F}_q) \}.$$

In order to simplify the notation, we shall use  $\nu (\bar{f}), \sigma (\bar{f})$ instead of $ \nu(\bar{f},D),  \sigma (\bar{f},D)$, respectively. The dependence of a particular set $D$, it will be clear from the constext.

With all this, we are able to establish  Igusa's stationary phase formula for $\pi-$adic integrals:
\bigskip
\subheading{ Igusa's Stationary Phase Formula} ([10], p. 17)

$$\int_{D}|f(x)|_{K}^{s}\mid  dx \mid = \nu(\bar{f}) +\sigma(\bar{f}) \frac{(1-q^{-1})q^{-s}}{(1-
q^{-1}q^{-s})} + \sum_{P \in S(f,D)} q^{-n -
e_{P}s}\int_{\Cal{O}_K^{n}} |f_{P}(x)|_K^{s}\mid dx \mid, \tag{2.3}$$

\noindent where $Re(s) >0$. Formula (2.3) is obtained in the following form. Suppose that 
  $\overline{D}= \{\overline{P_1},..,\overline{P_N}\}$
and let $P_i$ be the lifting of $\overline{P_i}$. Then the set $D$ is the disjoint union
$\bigcup_{P}D_P$, where $\overline{P}= (\overline{y_1},..,\overline{y_n}) \in \overline{D}$ and
$D_P$ is defined as
$$D_P = \{(x_1,...,x_n) \in D | \,\,\ x_i= y_i + \pi z_i, \,\,  z_i \in \Cal{O}_K, \,\ i=1,2,..,n  \}.$$
Thus

$$\int_{D}|f(x)|_{K}^{s}\mid dx \mid = \sum_{ \bar{P}} \int_{D_P} |f(x)|_{K}^{s}\mid dx \mid
 = \sum_{ \bar{P}}  q^{-n-e_{P}s} 
\int_{\Cal{O}_{K}^{n}}|f_{P}(x)|^{s}_{K}\mid dx \mid.$$

\noindent  The integrals corresponding to the $P$'s for which
 $\overline{P} \notin V_{\bar{f}}
 (\Bbb{F}_q)$
 are easily computable. The integrals corresponding to the $P$'s for which 
$\overline{P}$ is a non-singular point of $V_{\bar{f}}(\Bbb{F}_q)$ are computed using the 
implicit function theorem (cf. [10], p. 177). 
 
By iterating the stationary phase formula, we obtain the following  expansion for $Z(f,s)$ (cf. [10], p. 178):

$$Z(f,s) =\sum_{k \geq 0}q^{-kn}\left( 
\sum_{P_1,P_2,...P_k}\nu(\overline{f}_{P_1,...P_k})q^{-E(P_1,...P_k)s}\right)$$
$$ +\frac{(1-q^{-1})q^{-s}}{(1-q^{-1}q^{-s})} \sum_{k \geq 0}q^{-kn}\left( 
\sum_{P_1,P_2,...P_k}\sigma(\overline{f}_{P_1,...P_k})q^{-E(P_1,...P_k)s}\right)\tag{2.4}$$

\noindent where $E(P_1,..,P_k):= 
e_{P_1}+...+e_{P_1,...P_k}$. Expansion (2.4) converges 
absolutely on the semiplane $Re(s) > 0$.

Now, we  summarize some ideas on N\'{e}ron $\pi-$desingularization (see [13], sect. 17, 18) to be used  in the next
sections. Let $f(x) \in \Cal{O}_K[x] $ be a polynomial,  and  $P  \in V_{f}(\Cal{O}_K)$.  N\'{e}ron  introduced the following  measure of  singularity at $P$:

$$l(f,P) := \text{Inf}_{i} \left( v( \frac{\partial f}{\partial 
x_i}(P))\right).$$

\noindent The Jacobian criterion implies that  $ P$ is a smooth point of $V_{f}(K)$
 if and only if $l(f,P)$ is finite. $\overline{P}$ is a smooth
point of $V_{\bar{f}}(\Bbb{F}_q) $ if and only if $l(f,P)=0$.
In this paper, we introduce the following measure of singularity
 at an integer point $P$,  satisfying $\overline{P} \in V_{\bar{f}}(\Bbb{F}_q)$.   

\definition{Definition 2.1} Let $f(x)\in \Cal{O}_K[x]$ be
 a polynomial and $ P \in \Cal{O}_K^{n}$ a point, such that $ P \notin  Sing_{f}(\Cal{O}_K)$ 
 and $ \overline{P}  \in V_{\bar{f}}(\Bbb{F}_q)$. We define

$$L(f,P) := \text{Inf}
 \left( v(f(P)), v( \frac{\partial f}{\partial x_1}(P)),...,
 v( \frac{\partial f}{\partial x_n}(P))\right).$$
\enddefinition

\noindent Let us observe that   $L(f,P) =0$ if an only if 
$$  \overline{f(x)} = \alpha_0 +
\sum_{j} \alpha_j ( x_j - \overline{a_j}) + ( \text{ degree }
 \geq 2),$$
\noindent  where  $P = (a_1,..,a_n)$,  $\alpha_0 \in
  \Bbb{F}_q^{*}$ or  $\alpha_j \in  \Bbb{F}_q^{*}$ for some $j=1,2,..,n$.
We also observe  that $l(f,p) =L(f,p)$ if $P \in V_f(\Cal{O}_K)$.  
The integer $L(f,P)$ has similar properties to those of $l(f,P)$. This integer appears naturally associated to  Igusa's stationary phase, as we shall see later on.

We denote by $A_r$, $r= (r_1,..,r_n) \in (\Bbb{N}\setminus \{0\})^{n} $, the set
 $$A_r :=\{ x \in \Cal{O}_K^{n} \,\, \mid \,\,\, v(x_i) \geq r_i, \,\, i=1,..,n\}.$$
\noindent From a geometrical point of view, $A_r$ is a polydisc in $ \Cal{O}_K^{n}$ centered 
at the origin. The complement of $A_r$ in $\Cal{O}_K^{n}$ is denoted as $A_r^{c}$.

The following proposition is a simple reformulation  of proposition 17 in sect. 17
 of [13]. However for our convenience, we prove it below.

  \proclaim {Proposition 2.2}(N\'{e}ron, [13], sect. 17, prop. 17) Let $f(x)\in \Cal{O}_K[x]$
   be a  polynomial,  $P\in \Cal{O}_K^{n}$  an absolutely algebraically isolated singularity of
the hypersurface $V_f$,
   and let 
    $D \subseteq \Cal{O}_K^{n}$ be a subset such that $D \bigcap (P + A_r) = \emptyset$, for some $r \in (\Bbb{N} \smallsetminus \{0\})^{n}$. Then

    $$ L(f,Q) \leq C(f,D), \,\ \text{for every } \,\, Q \in D,$$
     
\noindent where the constant $C(f,D)$ depends only on $f$ and $D$.
\endproclaim
 \demo{Proof}
Without loss of generality, we may suppose that the point $P$ is the origin of $K^{n}$. The hypothesis that the origin of $K^{n}$ is an  absolutely algebraically isolated singularity and the Hilbert Nullstellensatz imply that
$$
\pi ^{m_i} x_i^{t_i}= A_{i,0}(x)f(x) + \sum_{j=1}^{n} A_{i,j}(x) \frac{\partial f}{\partial x_j}(x),\tag{2.5}$$
\noindent for  some $m_i$, $t_i \in  \Bbb{N}$ and some polynomials  $A(x)_{i,j} \in \Cal{O}_K[x]$, for each $i =0,1,2,..,n$. Now, let  $Q= (q_1,..,q_n)$ be a point of  $D$. Since 
$D \bigcap A_r = \emptyset$, there exists a coordinate $j_0$ such that
$v(q_{j_0}) < r_{j_0}$. From (2.5), with $x=Q$ and  $i=j_0$, we obtain

$$m_{j_0} +t_{j_0}r_{j_0} \geq m_{j_0}+ t_{j_0} v(q_{j_0}) \geq L(Q,f).$$
Thus, it is sufficient to take $C(f,D) \geq Max_i\{r_i + t_i m_i\}$.\quad\qed

\enddemo        
  The following result is a  generalization  of  proposition 18 (cf. [13],
   sect. 18) of N\'{e}ron.  

  \proclaim{Proposition 2.3}(N\'{e}ron, [13], sect. 18, prop. 18)
   Let $f(x)\in \Cal{O}_K[x]$
   be a  polynomial and $P \in \Cal{O}_K^{n}$ a point such that
    $P \notin Sing_{f}(\Cal{O}_K)$,  and   
    $\overline{P} \in Sing_{\bar{f}}(\Bbb{F}_q)$. Then there exists a minimal non-negative  integer
    $\mu(f,P)$, such that the polynomial

    $$   f_{P}(x)=\pi^{-e_{\mu ,P}}f(P + \pi^{\mu} x),$$

  \noindent where $e_{\mu ,P}$ is the minimum order of $\pi$ in the
  coefficients of $ f(P + \pi^{\mu} x)$, satisfies

  $$ \overline{f_{P}(x)}  =
  \alpha_0, \,\,\,\,\,\,\ \alpha_0 \in  \Bbb{F}_q^{*} \,\,\,\,\  \text{or} \,\,\,\,\ 
  \overline{f_{P}(x)}  =  \sum \alpha_i x_i, \,\,\,\,\,\  \alpha_i \in
   \Bbb{F}_q^{*}, \,\,\, \text{for some} \,\  i=1,2,..,n. $$

\noindent Furthermore, $\mu(f,P) \leq L(f,P)+2 $.
\endproclaim
\demo{Proof}
Let $P = (b_1,..,b_n) \in \Cal{O}_K^{n}$ be a point such that $\overline{P} \in
Sing_{\bar{f}}(\Bbb{F}_q)$.  Since $P \not \in Sin_f(\Cal{O}_K) $,
 we have
$$f(x) = \alpha_0 + \sum_i \alpha_i (x_i -b_i) + (\text{degree} \geq 2),$$
\noindent where $\alpha_i \equiv 0$  mod $\pi$,  $i=0,1,..,n$. Thus

$$f(P + \pi x) =\pi \left(\alpha_0' + \sum \alpha_i x_i +
\pi (\text{degree} \geq  2) \right).$$

\noindent We consider two cases according to $\alpha_0 '
  \not\equiv
    0$  mod $\pi$ or not.

 \noindent {\it Case 1}$ \,\,\,\, \text{(} \alpha_0 ' \not\equiv 0 $ mod $\pi$).

  \noindent In this case, we have
 $$f(P + \pi x)= \pi f_{P}(x),$$
 \noindent where

 $$f_P(x) = \alpha_0' + \sum \alpha_i  x_i  +
\pi (\text{degree} \geq 2).$$

\noindent Therefore, $\overline{f_{P}(x)}  = \overline{\alpha_0 '}\in  \Bbb{F}_q^{*}$,
and $\mu(f,P) =1 \leq  L(f,P)$.

\noindent  {\it  Case 2}    $\,\,\,\, \text{(} \alpha_0 ' \equiv  0$  mod $\pi$).

\noindent  In this case, we have

 $$f(P +\pi x) = \pi^{2} \left(\alpha_0'' + \sum_i \alpha_i'  x_i  +
 (\text{degree} \geq 2) \right)= \pi^{ e_{\mu, P} }f_P(x),$$
\noindent where $e_{\mu, P} \geq  2$. Thus

$$f(P) =  \pi^{e_{\mu, P}}f_P(0),$$
\noindent and
$$\frac{\partial f}{\partial x_i}(P) =
 \pi^{e_{\mu, P}-1}\frac{\partial f_P}{\partial x_i}(0), \,\,\, i=1,..,n.$$
 \noindent Whence $L(f_P,0) \leq L(f,P) -1$.
    Thus  after a finite number of dilatations, we obtain
    $L(f_{P_1,..,P_k},0) = 0$, ( where $P_2,..,P_k$ are equal to the origin of $K^{n}$), i.e.
    $$f_{P_1,...,P_k}(x) = \alpha_0 + \sum \alpha_i x_i +
     (\text {degree} \geq 2),$$
    \noindent where $\alpha_0 \not\equiv 0$ mod $\pi$ or $\alpha_i \not\equiv 0$
     mod $\pi$, for some $i$, $1 \leq i \leq n$.  If
     $\alpha_0 \not\equiv 0$  mod $\pi$, then after an additional dilatation, we
      obtain $\overline{ f_{P_1,...,P_k}(x) }=
    \overline{\alpha_0 } \in \Bbb{F}_q^{*}$, and $ \mu(f,P) \leq L(f,P)+1$.
     If $\alpha_0 \equiv  0$  mod $\pi$, then  after  an additional
 dilatation at the origin, we obtain
    
$$\overline{ f_{P_1,...,P_k}(x) }= \sum \overline{\alpha_i} x_i,$$
\noindent where $\overline{\alpha_i} \ne 0$ for some $i$, and
$\mu(f,P) \leq L(f,P) +2$.
\quad\qed
\enddemo

As a consequence of the two above results, we obtain the following lemma.

\proclaim{Lemma 2.4}
\noindent Let $f(x) \in \Cal{O}_K[x]$ be a polynomial such that the origin of $K^{n}$ is an absolutely algebraically isolated singularity. Let   $ A_r $ be  a polydisc with $r \in (\Bbb{N} \smallsetminus \{ 0 \})^{n}$. Then  there exists $\gamma= \gamma(f,r) \in \Bbb{N}$, such 
that the polynomial
$$f_Q(x) = \pi ^{-e _{Q, \gamma}}f(Q + \pi ^{\gamma}x),$$
satisfies the condition, $\overline {f_Q(x)}$ is a non-zero constant or a linear polynomial without constant term, for all $Q \in A_r ^{c}$. 
\endproclaim

\head
{\bf 3. Rationality of Igusa's local zeta functions of
semiquasihomogeneous polynomials}
\endhead
 
 In this section, we prove   the rationality of the Igusa
local  zeta function of  semiquasihomogeneous polynomials. 
 
\proclaim{Lemma 3.1}
 Let  $D \subseteq \Cal{O}_K^{n}$ be
the preimage under the canonical homomorphism $\Cal{O}_K \longrightarrow  \Cal{O}_K / \pi 
\Cal{O}_K $ of a subset $\overline{D} \subseteq \Bbb{F}_q^{n} $.  
Let
 $f(x) \in \Cal{O}_K[x] $ be a
 polynomial such that the origin of $K ^{n}$ is an absolutely algebraically isolated singularity of $V_f(K)$. If $D \bigcap A_r = \emptyset$, for some  $r \in (\Bbb{N} \smallsetminus \{0\})^{n}$, then 
the integral $Z(f,D,s) =\int_{D}|f |_{K}^{s}\mid dx \mid$ is a rational function
of  $q^{-s}$.  More precisely,       

$$ Z(f, D, s) = \frac{L(q^{-s}, D)}{1-q^{-1}q^{-s}}. \tag{3.1}$$

\noindent Furthermore, the polynomial $L (q^{-s},  D)$ can be effectively
 computed. 
\endproclaim
\demo{Proof} 
Applying the stationary phase formula  $m+1-$times, we obtain
\bigskip
\bigskip
$$Z(f,D,s) =\sum_{k= 0}^{m}q^{-kn}\left( 
\sum_{P_1,P_2,...P_k}\nu(\overline{f}_{P_1,...P_k})q^{-E(P_1,...P_k)s}\right)+$$
$$ \frac{(1-q^{-1})q^{-s}}{(1-q^{-1}q^{-s}} \sum_{k = 0}^{m}q^{-kn}\left( 
\sum_{P_1,P_2,...P_k}\sigma(\overline{f}_{P_1,...P_k})q^{-E(P_1,...P_k)s}\right)$$
$$+\sum_{P_1,..,P_m}\,  \sum_{P_{m+1} \in S(f_{P_1,..,P_m})}q^{-(m+1)n-E(P_1,...P_m)s}
  \int_{\Cal{O}_K^{n} }|f_{P_1,..,P_m}(x) |_{K}^{s}\mid dx \mid. \tag{3.2}$$

\noindent  On the other hand, we  have
$$f(P_1 + P_2\pi +....+ P_{m} \pi^{m+1}+ \pi^{m+2} x) =
\pi^{E(P_1,...,P_{m})}
   f_{P_1,....,P_{m}}(x). \tag{3.3}$$

Since $P_1 + P_2\pi +....+ P_{m} \pi^{m+1} \in A_r ^{c}$,  lemma 2.4 and (3.3) imply that 
 if $m+2 \geq \gamma(f, r)$, then 
the $\Bbb{F}_q-$hypersurface defined by
 $\overline{f_{P_1,....,P_{m}}(x)}$ is smooth or empty, whence the corresponding integral in (3.2) can be computed using the stationary phase formula.
Therefore  the integral
$ Z(f,D,s)$ can be computed applying $\gamma(f, r)$ times
the stationary phase formula. \quad\qed
\enddemo

\proclaim{Lemma 3.2}
Let
 $f(x) \in \Cal{O}_K[x] $ be a
 polynomial such that the origin of $K ^{n}$ is an absolutely algebraically isolated singularity of $V_f(K)$. 
Then the integral  $ Z(f,A_r^{c},s) = \int_{A_r^{c}}
 |f|^{s}\mid dx \mid$
   is a rational function of $q^{-s}$. More 
precisely,

$$Z(f,A_r^{c},s) = \frac {L(q^{-s})}
{
1-q^{-1}q^{-s}
}
.\tag{3.4}$$

\noindent  Furthermore, the 
polynomial $L(q^{-s})$ can be computed effectively.
\endproclaim

\demo{ Proof }
 We introduce a family $\Cal{L}$
of sets defined as follows.  For each   subset $B$ of $\{1,...,n\}$
 and each $ a = (a_1,...,a_n) \in \Bbb{N}^{n}$
satisfying $ 0 \leq a_i < r_i$ if $i \in B$, we define

$$ D(B,a):= \{ x \in A_r^{c} \mid \,\ v(x_i)=a_i \,\, \text{if} \,\  i \in B \},  \,\,\ 
\text{if} 
\,\, 
B 
\neq 
\emptyset,$$

$$ D(B,a):= \emptyset, \,\,\, \text{if} \,\, B= \emptyset.\tag{3.5}$$

\noindent The family $\Cal{L}$ is closed under  intersections, and its
 union is $A_r^{c}$.
 We denote by  $J$ the set of indices  $ \{(B,a)\}$ and by
 $ \Cal{P}(J)_i$  the family of subsets of $J$ with $i$ elements.

Whence
$$Z(f, A_r^{c},s) = \sum_{i=1}^{\text{Card} \{J\}} (-1)^{i-1}\sum_{T \in \Cal{P}(J)_i}
\int_{ D(T) }
|f|^{s}\mid dx \mid ,\tag{3.6}$$
where $D(T):= \bigcap_{(B,a) \in T}D(B,a)$.
From (3.6) and the fact that the family $\Cal{L}$ is closed under intersections, it  follows that  in order to prove the theorem, it is sufficient to prove  that
any integral of type
$$\int_{D(B,a)}|f|^{s} \mid dx \mid,  \,\,\ B \neq \emptyset \tag{3.7}$$
 is a rational function
of the form $L(q^{-s}, D(B,a))/(1-q^{-1}q^{-s})$, where the numerator
 polynomial
is effectively computable. For that, we make the following change of  variables in (3.7),
 $$x=  \psi_{(B,a)}(y), \,\, \text{where} \,\
 x_i=
 \cases
  \pi^{a_i}y_i \,\ \text{if} \,\  i \in B\\
y_i \,\,\ \text{if} \,\  i \notin B.\\ 
 \endcases
\tag{3.8}
 $$
  we obtain
 $$\int_{D(B,a)}|f|^{s} \mid dx \mid = q^{-e_{(B,a)}s- d_{(B,a)} }
 \int_{D^{'}(B,a)}|f_B|^{s} \mid dy \mid, \tag{3.9}$$
 \noindent where $d_{(B,a)} = \sum_{i \in B} a_i $,  $f_B(y)$ is the dilatation of $f$ induced by (3.8)
 and
 $$D^{'}(B,a) = \prod_{i=1}^{n}R_i,$$
 where $R_i =\Cal{O}_K$ if $i \notin B$
 and  $R_i =\Cal{O}_K^{*}$ if $i \in B$.

On the other hand, $\phi(y)$ defines a $K-$isomorphism of
$K^{n} \longrightarrow K^{n}$, thus the $K-$singular locus of $V_f$
is mapped bijectively on the $K-$singular locus of $V_{f_B}$. Therefore, the polynomial
 $f_B$   and the set $D^{'}(B,a)$, $B \neq \emptyset$, satisfy the conditions of lemma 3.1.
 Thus
 the integral  $\int_{D^{'}(B,a)}|f_B|^{s}\mid dx \mid$ is a rational function of $q^{-s}$ and its numerator can be computed effectively.\quad\qed
\enddemo

 \proclaim{Proposition 3.3} 
Let $F(x)= f(x) + \pi^{m}  g(x) \in \Cal{O}_K[x]$ be a semiquasihomogeneous
  polynomial ,  $f(x)$ is its quasihomogeneous part.  
 Let  $D \subseteq \Cal{O}_K^{n}$ be the preimage under the canonical
  homomorphism $\Cal{O}_K \longrightarrow  \Cal{O}_K / \pi \Cal{O}_K $ of a subset $\overline{D} \subseteq \Bbb{F}_q^{n} $, and $A_r$ a polydisc, 
such that $D \bigcap A_r = \emptyset$,   $r \in (\Bbb{N} \smallsetminus \{0\})^{n}$.
 There exists  $\alpha(f,D)$, effectively computable,  such that if
  $m \geq \alpha(f,D)$ then  
 
 $$Z(F,D,s) = Z(f,D,s). $$
\endproclaim

\demo{Proof} 
By virtue of  lemma 2.4, there exists a $\gamma _0$ such that the reduction modulo
$\pi$ of the polynomial
$$F_{P_1,...,P_n}(x) = \pi^{ -E_{F}(P_1,..,P_{n} )}F( P_1 + P_2\pi +....+P_n \pi^{n} +\pi^{n+1}x)\tag{3.10}$$
is a non-zero constant or a linear polynomial for every $n \geq \gamma _0$ and any
$P_1 + P_2\pi +....+P_n \pi^{n} \in A_r^{c}$. In addition, we also have that
$$
F( P_1 + P_2\pi +....+P_n \pi^{n} +\pi^{n+1}x) =\pi^{ E_f(P_1,..,P_{n} )}f_{P_1,...,P_n}(x)$$
$$ +  \pi^{ E_g(P_1,..,P_{n} )+ m}g_{P_1,...,P_n}(x).\
\tag{3.11}$$

We choose 
$$ \alpha(f,D) := Max \left \{ E_f(P_1,..,P_{\gamma _0}) \right\},\tag{3.12}
$$
where  $P_{\gamma _0}$ runs through the $\gamma_0$ generation of decendants
 of $S(F,D)$.
Now,  If $ m> \alpha(f,D)$, we have
 $$\overline{F_{P_1, P_2,..,P_k}(x)}=\overline{f_{P_1, P_2,..,P_k}(x)},$$
 $$E_{P_1, P_2,..,P_k }(F) = E_{P_1, P_2,..,P_k }(f),  \,\,\,\,\, 1 \leq k \leq \gamma_0.\tag{3.13} $$
Expanding $Z(F,D,s)$ and $Z(f,D,s)$ as in (3.2) and using (3.10) and (3.13), we conclude that
 $Z(F,D,s) =Z(f,D,s).$
 
 \quad\qed
 \enddemo
 \proclaim{Lemma 3.4}
 Let $F(x)= f(x) + \pi^{m}  g(x) \in \Cal{O}_K[x]$ be a semiquasihomogeneous
  polynomial,  $f(x)$ is its quasihomogeneous part.  
 Let $ A_r $, be a polydisc, with 
 $r \in (\Bbb{N} \smallsetminus \{0\})^{n}$. 
    There exists  $\alpha(f,r)$, effectively computable,  such that if
  $m \geq \alpha(f,r)$ then  
 
 $$Z(F,A_r^{c},s) = Z(f,A_r^{c},s). $$
\endproclaim 

\demo{Proof} By (3.6), it is sufficient to prove the lemma for the integrals of type (3.7).
We choose $\alpha (f,r)$ satisfying 
$$\alpha (f,r) \geqq Max_{(B,a)} \{  e_{f,{(B,a)}}\}.$$
With the above condition, we have that (see (3.9))it is sufficient to prove that
$$\int_{D(B,a)'}|F_B|^{s}\mid dx \mid = \int_{D(B,a)'}|f_B|^{s}\mid dx \mid.\tag{3.14}$$
The result follows from (3.14) and proposition 3.3. Finally, we observe that $\alpha (f,r)$ is given by
$$
\alpha (f,r) = Max_{(B,a)} \{  e_{f,{(B,a)}}+ \alpha(f, D^{'}(B,a))\}.
$$ \quad\qed
 \enddemo

\proclaim{Theorem 3.5} 
Let $F(x) \in K[x]$ be an a semiquasihomogeneous
 polynomial whose quasihomogeneous part $f(x)$
has weight $d$
and exponents $ \alpha_1,..,\alpha_n$. Then   Igusa's local zeta function of $F(x)$ is a rational function of $q^{-s}$. More
precisely,

$$Z(F,s) = \frac{L(q^{-s})}{ (1-q^{-1}q^{-s})(1-q^{-|\alpha|}q^{-ds})}
.\tag{3.15}$$

\noindent where $\alpha= (\alpha_1,...,\alpha_n)$.  Furthermore, the 
polynomial $L(q^{-s})$ can be computed effectively.
\endproclaim

\demo{Proof}
By decomposing  $\Cal{O}_K^{n}$
as the disjoint union of $A_{\alpha}$ and $A_{\alpha}^{c}$ and using the fact that $F(x)$ is
a semiquasihomogeneous polynomial, we obtain

$$Z(F,s) =  \int_{A_{\alpha}}|F|_K^{s}\mid dx \mid + \int_{A_{\alpha}^{c}}|F|_K^{s}\mid dx \mid =  
q^{-|\alpha| -ds} \int_{ \Cal{O}_K^{n}}  |F_1|_K^{s}\mid dx \mid +
 \int_{A_{\alpha}^{c}}|F|_K^{s}\mid dx \mid, \tag{3.16}$$

 \noindent where $F_1(x) =f(x) +\pi ^{m_1}H_1(x)$, with $m_1 \geq 1$.
Now, Iterating  formula (3.16) $m-$times, we obtain:
$$Z(F,s) = Z(F, A_{\alpha}^{c},s) +
 \sum_{k=1}^{m} q^{k(-|\alpha| -ds)} Z(F_k,A_{\alpha}^{c},s)
  + q^{(m+1)(-|\alpha| -ds)}Z(F_{m+1},s), \tag{3.17}$$
where $F_k(x) = f(x) +\pi^{m_k}H_k(x)$, with $ m_k \longrightarrow \infty$.
By lemma 3.4, there exists $\gamma_0 = \gamma_0(f, \alpha)$, effectively computable,  such that $Z(F_k,A_{\alpha}^{c},s)= Z(f,A_{\alpha}^{c},s)$ if $ m \geq
\gamma_0$. Thus from (3.17), we have
$$
Z(F,s) = Z(F,A_{\alpha}^{c},s) + \sum_{k=1}^{m_0 -1}q^{k(-|\alpha| -ds)} Z(F_k,A_{\alpha}^{c},s) +
$$
$$
Z(f,A_{\alpha}^{c},s)q^{(m_0+1)(-|\alpha| -ds)}\frac{1}{1-q^{-|\alpha| -ds}}.
\tag{3.18}
$$
By lemma 3.2 the integrals  $Z(f,A_{\alpha}^{c},s)$  and $Z(F_k,A_{\alpha}^{c},s)$  are rational functions of $q^{-s}$ of the
 form $\frac{L(q^{-s} )}{1- q^{-1}q^{-s}}$, where the polynomial
  numerator can be  effectively computed.\quad\qed
\enddemo

  We observe  that if $f$ is a quasihomogeneous polynomial,  its local zeta function is given by $Z(f,s)  = \frac{ Z(f,A_{\alpha}^{c},s)}{1-q^{-|\alpha| -ds}}$. The integral  $Z(f,A_{\alpha}^{c},s)$ can be computed using lemma 3.2. 

As a consequence of theorem 3.5, we obtain the following three corollaries.

\proclaim{Corollary 3.6}
Let $K$ be a global field and  let $F(x) \in K[x]$ be a semiquasihomogeneous
 polynomial whose quasihomogeneous part $f(x)$  has  weigth $d$
and exponents $ \alpha_1,..,\alpha_n$. 
Then for every  non-archimedean valuation $v$ of $K$,  Igusa's local zeta function of $F(x)$ on the 
completion $K_v$ of $K$ is a rational function of form (3.15). If $K$ is a number field
and  $F(x)$ is non-degenerate for
 its Newton's diagram, then  the real parts of the poles of the zeta function $Z(F,s)$  are roots of the Bernstein polynomial of $F(x)$.
\endproclaim
\demo{Proof}
Since $F(x)$ has an  absolutely algebraically isolated singularity at the origin of $K^{n}$,
 the Hilbert Nullstellensatz
 implies that for  all valuations $v$ of $K$, the origin of $K_v^{n} $ is  an  absolutely algebraically isolated singularity of $V_{f}(K_v)$. Thus by the proof of theorem 3.5, the denominator of the
 the local zeta function $Z(f,s)$ on $K_v$ is equal to $(1-q^{-1}q^{-s})(1-q^{-|\alpha|}q^{-d s})$. Thus the real parts of the  poles of $Z(F,s)$ are among the values $-1$, $- \mid \alpha \mid/d$.  If  $K$ is   a number field,  
and  $F(x)$ is non-degenerate for its Newton's diagram,
 theorem   C.2.2.3 of [2] implies that $-1$ and $ -\mid \alpha \mid/d$ are roots of the Bernstein polynomial of $F(x)$. \quad\qed
 \enddemo
The following corollary gives  a bound for the number of solutions  of a conguence attached to
a semiquasihomogeneous polynomial with coefficients in a ring of integers of a global field.
This corollary follows directly from the relation existing between the Igusa local zeta function and the Poincar\'{e} series $P(t)$ (see introduction) and corollary 3.6.

\proclaim{Corollary 3.7}
\noindent Let $K$ be a global field and $\Cal{O}_K$ its ring of integers. Let $F(x) \in \Cal{O}_K[x]$ be a semiquasihomogeneous polynomial whose quasihomogeneous part $f(x)$ has weight $d$ and exponents $\alpha_1,...,\alpha_n$. Then for every non-archimedean 
valuation $v$ of $K$, the number of solutions $N_j(F,v)$ of the congruence
$$ F(x)
 \equiv 
0 
\,\ 
\text{mod} 
\,\  \pi^{j}\Cal{O}_{K_{v}},$$
where  $\Cal{O}_{K_{v}}$  is the ring of integers of the completion $K_v$, satisfies 
$$
\text{l\'\i m sup}_{j\to \infty}N_{j}(F,v)^{1/j} \le 
\cases
 q^{n- \mid \alpha \mid  / d }\,\,\,\,\ \text{if}  \,\,\ \mid \alpha \mid  / d \le 1, \,\,\ \\
q^{n- 1 } \,\,\,\,\,\,\,\,\,\,\,\  \text{if} \,\,\,\, \mid \alpha \mid  / d > 1.\\
\endcases
$$
\endproclaim
 \medskip
 The following corollary follows directly from theorem 3.5. We use the notation established in the introduction.

\proclaim{Corollary 3.8}

\noindent Let   $f(x,y) \in K[x,y]$ be an absolutely  analytically irreducible polynomial, such that the origin of $K^{2}$ is an absolutely algebraically isolated singularity of
 the plane
curve   $V_{f_d}$.  Then  the Igusa local zeta function $Z(f,s)$ is a rational function of
$q^{-s}$ of  form (3.15).
\endproclaim

\example {Example 3.9} 
In this example  we compute the local zeta function for a  polynomial of type  $f(x,y) = \alpha x^{n} + \beta y^{m}$, $\alpha, \beta \in \Cal{O}_K$, where $n,m >1$ are relatively prime. 
   Suppose that  the characteristic of $K $ does not divide both  $n,m$. Furthermore, without loss of generality, we may suppose that
$\alpha
 \in \Cal{O}_K^{*}$. In the case of charactersitic zero,
the Poincar\'{e}
  series $P(t)$ associated to this type of  polynomials were explicitly
   computed by Goldman (cf. [6], thm. 1).

Using the observation made after the proof of theorem 3.5, we have
   
$$Z(f,s) =\frac{1}{1-q^{-(n+m) -mns}}\int_{A^{c}} |f|^{s}_{K}\mid dx dy \mid ,\tag{3.19}$$

\noindent where $A = \{(x,y) \in \Cal{O}_{K}^{2} \,\,\ | \,\,\ v(x) \geq m, \,\,\,\, v(y) \geq n \}$.The complement $A^{c}$ of $A$ is the disjoint union of the following three sets:

$$D_1 = \{(x,y) \in \Cal{O}_{K}^{2} \,\,\ | \,\,\ v(x) < m, \,\,\,\, v(y) \geq n \},$$
$$D_2 = \{(x,y) \in \Cal{O}_{K}^{2} \,\,\ | \,\,\ v(x) < m, \,\,\,\, v(y) < n \}.$$
$$D_3 = \{(x,y) \in \Cal{O}_{K}^{2} \,\,\ | \,\,\ v(x) \geq m, \,\,\,\, v(y) < n \}.$$
Thus from (3.19), we obtain
$$Z(f,s) =\frac{1}{1-q^{-(n+m) -mns}} \left \{ Z(f,D_1,s) +
 Z(f,D_2,s) + Z(f,D_3,s) \right \}
.\tag{3.20}$$

Next, we compute the integrals $Z(f,D_1,s)$, $Z(f,D_2,s)$, $Z(f,D_3,s)$.

\noindent {\bf Computation of $Z(f,D_1,s)$}

First, we observe that
$$\mid f(x,y) \mid = \mid  \alpha x^{n} + \beta y^{m} \mid = \mid x^{n} \mid, \,\
\,\,\,\  x,y \in D_1.$$
\noindent  Therefore

$Z(f,D_1,s) = \int_{D_1}|f|^{s}_{K}\mid dx dy \mid =
 \sum_{k=0}^{m-1} \int_{ \{(x,y) \in D_1\,\ | \,\ v(x) = k, \, v(y) \geq n \}}
 |x|^{ns}_{K}\mid dx dy \mid.$

\noindent Thus
$$Z(f,D_1,s)= (1-q^{-1})q^{-n} \sum_{k=0}^{m-1} q^{-kns-k}\tag{3.21}$$

\noindent {\bf Computation of $ Z(f,D_2,s) $}

We set 
 $ L(i,j):= jm -in +v(\beta)$.
The set $D_2$ can be decompossed as  the union of three disjoint subsets $D_{2,1},D_{2,2},D_{2,3}$, as 
follows :

$$D_{2,1} :=\{(x,y) \in D_2 \,\ | L( v(x),v(y)) >0 \},$$
$$D_{2,2} :=\{(x,y) \in D_2 \,\ | L( v(x),v(y)) <0 \},$$
$$D_{2,3} :=\{(x,y) \in D_2 \,\ | L( v(x),v(y)) =0 \}.$$
\noindent Thus  $Z(f,D_2,s) =  Z(f,D_{2,1},s)+Z(f,D_{2,2},s)+
Z(f,D_{2,3},s)$, where

$$Z(f,D_{2,1},s) =(1-q^{-1})^{2}\sum_{i,j}q^{-i-j  - nis},\tag{3.22}$$
\noindent where $i,j$ satisfy $L(i,j) >0$ and $ 0 \leqq i <m$, $0 \leqq j <n$,

$$Z(f,D_{2,2},s) ={(1-q^{-1})}^{2}\sum_{i,j}q^{-i-j -(v(\beta) + mj)s},\tag{3.23}$$
\noindent where $i,j$ satisfy $L(i,j) <0$ and $ 0 \leqq i <m$, $0 \leqq j <n$, and

$$Z(f,D_{2,3},s) =\sum_{i,j}q^{-i-j -nis}\int_{\Cal{O}_{K}^{\times \, 2}}|\alpha x^{n} + \mu  y^{m}|^{s}_{K}\mid dx dy \mid,\tag{3.24}$$
\noindent where $\beta= \pi ^{v(\beta)}\mu$, $\mu \in \Cal{O}_K^{*}$, $i,j$ satisfy $L(i,j) =0$ and $ 0 \leqq i <m$, $0 \leqq j <n$.
 Using the stationary phase formula, we compute the integral in the right side  of  (3.24), thus 

$$Z(f,D_{2,3},s) = \left (\nu(\overline{f})  + 
\frac{\sigma(\overline{f})(1-q^{-1})q^{-s}}{1-q^{-1-s}}\right)
 \sum_{i,j} q^{-i-j - nis} .$$

\noindent We denote by $[x]$ the integer part of  a  real number $x $.We set $v(\beta) = gn + r$,
$ 0 \leq r < n$.

\noindent{\bf Computation of $ Z(f,D_3,s) $}

We set 
$$D_{3,1} := \{(x,y) \in \Cal{O}_K^{2} \,\,\ | \,\,\ v(x) \geq m + [\frac{v(\beta)}{n}]+ r, \,\,\,\, v(y) < n \},$$

$$D_{3,2} := \{(x,y) \in \Cal{O}_K^{2} \,\,\ | \,\,\ m \leq v(x) \leq m + 
[\frac{v(\beta)}{n}]+r -1, \,\,\,\, v(y) < n \}.$$

Then $D_3 = A_{3,1} \bigcup D_{3,2}$ (disjoint union), and $Z(f,D_3,s) =
Z(f,D_{3,1},s) + Z(f,D_{3,2},s)$. 
To compute $Z(f,D_{3,1},s)$, we observe that
$$| f(x,y)| = |\alpha x^{n} + \beta y^{m} |= | \beta y^{m}
 |\,\,\,\, x,y \in D_{3,1}.$$
\noindent Thus

$$Z(f,D_{3,1},s) = \int_{D_{3,1}}| f(x,y)|_K^{s}\mid dx  dy \mid = 
(1-q^{-1})q^{-(m+ [\frac{v(\beta)}{n}]+r)}\sum_{k=0}^{n-1}q^{-(v(\beta) +mk)s-k}.
\tag{3.25}$$

The set $D_{3,2}$ can be decompossed as  the union of three disjoint subsets  $D_{3,2,1}$,  $D_{3,2,2}$, $D_{3,2,3}$, as 
follows :

$$D_{3,2,1} :=\{(x,y) \in D_{3,2} \,\ | L( v(x),v(y)) >0 \},$$
$$D_{3,2,2} :=\{(x,y) \in D_{3,2} \,\ | L( v(x),v(y)) <0 \},$$
$$D_{3, 2,3} :=\{(x,y) \in D_{3,2} \,\ | L( v(x),v(y)) =0 \}.$$
\noindent Thus  $Z(f,D_{3,2},s) =  Z(f,D_{3,2,1},s)+Z(f,D_{3,2,2},s)+
Z(f,D_{3,2,3},s)$, and

$$Z(f,D_{3,2,1},s) =(1-q^{-1})^{2}\sum_{i,j}q^{-i-j  - nis},\tag{3.26}$$
\noindent where $i,j$ satisfy $L(i,j) >0$ and $ m \leqq i < m +[\frac{v(\beta)}{n}]+r $, 
$0 \leqq j <n$,

$$Z(f,D_{3,2,2},s) ={(1-q^{-1})}^{2}\sum_{i,j}q^{-i-j -(v(\beta) + mj)s},\tag{3.27}$$
\noindent where $i,j$ satisfy $L(i,j) <0$ and $ m \leqq i <m +[\frac{v(\beta)}{n}]+r $, $0 \leqq j <n$, and

$$Z(f,D_{3, 2,3},s) = \left (\nu(\overline{f})  + 
\frac{\sigma(\overline{f})(1-q^{-1})q^{-s}}{1-q^{-1-s}}\right)
 \sum_{i,j} q^{-i-j - nis},$$

where $i,j$ satisfy $L(i,j) = 0$ and $ m \leqq i <m +[\frac{v(\beta)}{n}]+r $, $0 \leqq j <n$.

\endexample

\example{Example 3.10}
A polynomial of the form $f(x) = \sum_{i}\alpha_i x_i^{n_i}$, $\alpha_i \in \Cal{O}_K$
is called a diagonal polynomial. We set $d:= l.c.m \{ n_i \}$, and
 $\alpha_i := d/n_i $, $i=1,..,n$. If the characteristic of $K$ does not divide any $n_i$, then the diagonal polynomials are quasihomogeneous
  polynomials with exponents   $\alpha_i := n_i / d$, $i=1,..,n$   and weight $d$.
 Thus  the local zeta function  of a  diagonal polynomial is
 a  rational function of  form (3.13). Wang and others have studied the Poincar\'{e} series $P(t)$ associated
to this class of polynomials  (cf [14], thm. 1). 
\endexample

  \subheading {Acknowledgements}

  The author wishes to thank to the following institutions for their support:
  Universidad Aut\'{o}noma de Bucaramanga, 
   Academia Colombiana  de Ciencias  Exactas, F\'{\i}sicas y Naturales
and  COLCIENCIAS. The author also
  thanks  IMPA  for their support and hospitality during the summer of 1997, when part of this work was done. The author also wishes to thank to   referee for his or her useful comments which lead to an improvement of  this work.
\refstyle{C}

\enddocument